\documentclass[12pt]{amsart}%
 \usepackage{graphicx}
\usepackage[mathscr]{eucal}
\makeindex
\newtheorem{Thm}{Theorem}[section]
\theoremstyle{definition}
\newtheorem{Theorem}[Thm]{Theorem}
\newtheorem{Lemma}[Thm]{Lemma}
\newtheorem{Corollary}[Thm]{Corollary}
\newtheorem{Proposition}[Thm]{Proposition}
\newtheorem{Definition}[Thm]{Definition}

\theoremstyle{remark}
\newtheorem{Remark}{Remark}
\font\sy=cmsy10

\font\ym=msbm10

\newcommand{\N}{{\text{\ym N}}}

\newcommand{\R}{\text{\ym R}}

\newcommand{\C}{\text{\ym C}}
\newcommand{\cB}{{\hbox{\sy B}}}

\newcommand{\sA}{\mathscr A}

\newcommand{\sD}{\mathscr D}

\newcommand{\sN}{\mathscr N}

\newcommand{\sR}{\mathscr R}

\newcommand{\Hom}{\hbox{\rm Hom}}

\title[analytic functionals]{On moment problems of analytic functionals of polynomial hypergroups}
\author{Shigeru Yamagami}
\begin{document}
% Haagerup positive definite functions on free groups.
% To Nobuko with our deep gratitudes. 
\maketitle
\begin{center}
Graduate School of Mathematics
\end{center}
\begin{center}
Nagoya University 
\end{center}
\begin{center} 
Nagoya, 464-8602, JAPAN
\end{center}

\subjclass{Primary 46J25; Secondary 46L99}

\keywords{moment problem, polynomial hypergroup, Haagerup function}
% , Stieltjes inversion formula}
% analytic linear functional
% 46L99 (operator algebra), 46J25(commutative Banach algebra), 43A99(abstract harmonic analysis) 
% Related to radial functions on free groups, we focus on certain polynomial hypergroups and
% work out spectral analysis with the help of the Stieltjes transform of their analytic functionals. 

\begin{abstract}
  Spectral analysis is worked out on moment problems related to linear functionals of certain polynomial hypergroups. 
\end{abstract}

% This work comes from free group harmonics: Cartier and Haagerup.
% Cartier did harmonic analysis in a form of hypergroup arizing from radial functions on free groups, whereas
% Haagerup introduced an interesting class of positive definite functions on free groups.
% Although there are lots of related results, a unified coverage of these seems to be still missing. 

\section*{Introduction}
We consider here a one-parameter family of polynomial hypergroups $\sA$,
which appeared implicitly in \cite{Ca} as an algebraic system of spherical functions of free groups
(cf.~also \cite[Chapter 3]{FP}),
% ($r = 1/2l$ in the case of radial functions of a free group)
and describe analytic linear functionals on $\sA$ as solutions to the associated moment problems.
Though our method is standard (computations of similar flavor can be found in \cite{Ak,Vo} for example),
it provides a new framework of spectral analysis together with algorithmic methods of computations.
% based on Stieltjes inversion formula together with generating functions,
% which is a common tool in dealing with problems of this kind,
% and computations of similar flavor can be found in 
% and investigate
% spectral properties of analytic linear functionals on $\sA$.
% % parametrized by $s>0$ which are obvious analogues of Haagerup's positive definite functions. 
% Our approach is standard and based on Stieltjes inversion formula together with generating functions. 
% Such a method is quite popular in dealing with recurrence relations 
% and computations of similar flavor can be found in \cite{Ak,Vo} for example.

As an illustrating example, the case of geometric series functionals is then explicitly worked out
in connection with Haagerup's positive definite functions on free groups and
we obtain a new series of spectral measures including the classical Kesten type measures and atomic masses,
which turns out to reflect some features of the original situations in free groups and
be useful in elucidating properties of free group C*-algebras.
As a concrete instance, Haagerup's positive definite functions are disintegrated with respect to these measures,
which includes the Plancherel formula of Cartier as a special case. 

% and calculate the representing Radon measures explicitly with C*-functionals selected out among them.
% Although we work in a completely commutative setting,
% the results turn out to reflect

The author is grateful to Yoshimichi Ueda for illuminating comments and useful information on the present subject.
% The author would also like to thank the referee for valuable comments, which inspired him to improve the appendix B. 

% and to the referee for drawing author's attention to \cite{BB}. 

\section{Polynomial Hypergroups}
  Given a real number $r \not= 1$, %$\alpha$ and $\beta \not= 0$,
  algebraic relations $tp_n(t) = r p_{n-1}(t) + (1-r) p_{n+1}(t)$ ($n \geq 1$)
  together with the initial condition $p_0(t) = 1$, $p_1(t) = t$
  determine polynomials $p_n(t)$ of real coefficients of degree $n$ inductively.
  Since $(p_n(t))$ constitutes a real (hermitian) basis in the polynomial *-algebra $\C[t]$,
  we can rephrase this fact in a slightly abstract fashion that a free vector space
  $\sum_{n \geq 0} \C h_n$ is made into a commutative *-algebra $\sA$ by the relations 
  $h_n^* = h_n$, $h_0 h_n = h_n$ ($n \geq 0$) and
  $h_1h_n = r h_{n-1} + (1-r) h_{n+1}$ ($n \geq 1$)
  so that $\sA$ is *-isomorphic to $\C[t]$ by the correspondence $h_n \leftrightarrow p_n(t)$.

  Notice that, for $0 \leq r < 1$,
  coefficients of $t^n$ (corresponding to $h_1^n$) relative to the basis $(p_0,p_1,\dots,p_n)$ are related to
  the reaching probabilities of a random walk on $\{0,1,2,\dots\}$: Here are several lower steps. 
  \begin{align*}
 %   t & = p_1\\
    t^2 &= rp_0 + (1-r) p_2,\\
    t^3 &= r(2-r) p_1 + (1-r)^2 p_3,\\
    t^4 &= r^2(2-r) p_0 + r(1-r)(3-2r) p_2 + (1-r)^3 p_4,\\
    t^5 &= r^2(5-6r+2r^2)p_1 + r(1-r)^2(4-3r)p_3 + (1-r)^4p_5.
  \end{align*}
  \begin{align*}
    p_0 &= 1,\\
    p_1 &= t,\\
    (1-r) p_2 &= t^2 - r,\\
    (1-r)^2 p_3 &= t^3 - r(2-r)t,\\
    (1-r)^3p_4 &= t^4 - r(3-2r)t^2 + r^2(1-r),\\
    (1-r)^4p_5 &= t^5 - r(4-3r)t^3 + r^2(1-r)(3-r). 
  \end{align*}
  As special values of $p_n$, we record
  \[
    p_n(\pm 1) = (\pm 1)^n,
    \quad
    p_{2n}(0) = \frac{(-r)^n}{(1-r)^n}.
  \]

% Given a real parameter $r \not= 1$ ($r = 1/2l$ in the case of radial functions of a free group),
% introduce polynomials $P_n(t)$ of degree $n$ by the recursive relation
% \[
%   tP_n(t) = rP_{n-1}(t) + (1-r) P_{n+1}(t) \quad (n \geq 1),
%   \quad
%   P_0(t) = 1,\ P_1(t) = t.
% \]
% As in the remark, we have a *-algebra $\sD = \sum_{n \geq 0} \C \delta_n$ so that
% $\delta_n \mapsto P_n(t)$ gives a *-isomorphism between $\sD$ and
% $\C[t]$ ($\{ \delta_n = \delta_n^*\}$ being a linear basis of $\sD$).
The multiplication $h_m h_n = h_n h_m$ is then recursively determined by the relation
$h_1 h_n = rh_{n-1} + (1-r)h_{n+1}$ ($n \geq 1$): 
% For example,
% $\delta_1^2 = r \delta_0 + (1-r)\delta_2$ is used to have 
% \begin{align*}
%   (1-r)\delta_2\delta_n &= -r \delta_n + \delta_1^2\delta_n
%   = - r\delta_n + \delta_1(r\delta_{n-1} + (1-r)\delta_{n+1})\\
%                         &= - r\delta_n + r(r\delta_{n-2} + (1-r) \delta_n) + (1-r)(r\delta_n + (1-r) \delta_{n+2})\\
%   &= r^2\delta_{n-2} + r(1-2r)\delta_n + (1-r)^2 \delta_{n+2} 
% \end{align*}
% and 
% \[
%   (1-r)^2 \delta_3 \delta_n = r^3\delta_{n-3} + r^2(1-2r)\delta_{n-1} + r(1-r)(1-2r) \delta_{n+1} + (1-r)^3\delta_{n+3}.
% \]
We can inductively show that, for $1 \leq m \leq n$, 
\begin{multline*}
  (1-r)^{m-1} h_m h_n\\
  = r^m h_{n-m} + r^{m-1}(1-2r) h_{n-m+2}
  + r^{m-2}(1-r)(1-2r) h_{n-m+4}\\
  + \cdots + r(1-r)^{m-2}(1-2r) h_{n+m-2} + (1-r)^m h_{n+m}. 
\end{multline*}
% The cancellation for the last term needs an exceptional care. 
Notice that the coefficients are positive (including $0$) exactly when $0 \leq r \leq 1/2$ (this parameter range being assumed in what follows)
and summed up to be $(1-r)^{m-1}$.
Thus $\sA$ is a discrete hypergroup % (simply a polynomial algebra when $r=0$)
in this range
% \footnote{Strictly speaking, hypergroups are for the range $0 < r \leq 1/2$.}
of the parameter and completed to a commutative Banach *-algebra relative to the norm
$\| \sum_{n \geq 0} a_n h_n\|_1 = \sum_{n \geq 0} |a_n|$.
See \cite{Ca} for the original form of $\sA$ and \cite{BH} for more information on (polynomial) hypergroups.
Let $A$ be the enveloping C*-algebra of this Banach *-algebra (cf.~\cite[Remark 6.1.1]{Mu}):
For a $\| \cdot\|_1$-bounded *-representation $\pi$ of $\sA$ on a Hilbert space, the operator norm $\| \pi(a)\|$ gives
a C*-seminorm on $a \in \sA$, whence $\| a\| = \sup_\pi \{ \| \pi(a)\|\}$ is a C*-seminorm as well.
Then $\sN = \{ a \in \sA; \| a\| = 0\}$ is a *-ideal of $\sA$ and 
$A$ is defined to be the C*-completion of the quotient *-algebra $\sA/\sN$ with respect to the induced C*-norm. 

{\small
\begin{Remark}\ 
  \begin{enumerate}
    \item
  For $r=1/2$, $p_n(t)$'s are Chebyshev polynomials.
\item
  For $r=0$, the polynomial basis $p_n(t) = t^n$ does not define a hypergroup but it is a limit of hypergroup bases.
\end{enumerate}
\end{Remark}
}

A linear functional $\varphi:\sA \to \C$ is naturally identified with a sequence $(\varphi_n)$ by
the relation $\varphi(h_n) = \varphi_n$, i.e.,
the algebraic dual $\sA^*$ of $\sA$ is identified with the space $\C^\N$ of sequences.
As a polynomial algebra, any non-zero homomorphism $\varphi: \sA \to \C$ satisfies $\varphi(h_0) = 1$ 
and is in one-to-one correspondence with a complex number $c = \varphi(h_1)$,
i.e., the evaluation of polynomials at $t = c$, which determines
the whole values $\varphi(h_n)$ by the recurrence relation
$c \varphi(h_n) = r \varphi(h_{n-1}) + (1-r) \varphi(h_{n+1})$ ($n \geq 1$). 
Among them, *-homomorphisms correspond to real $c$'s.

% By Gelfand theory (see, for example, \cite[Theorems 2.1.7 and 2.1.9]{Mu}), 
% *-homomorphisms which lift to $A$ are ones satisfying $|\varphi(a)| \leq \| a\|_1$
% ($a \in \sA$), i.e., $(\varphi(h_n))$ being bounded. % as an element in $(\ell^1)^*$.
Since any *-homomorphic functional $\varphi$ on $\sA$ lifts to $A$ if and only if $|\varphi(a)| \leq \| a\|_1$ ($a \in \sA$) and
any homomorphic functional on $A$ is *-preserving (see, for example, \cite[Theorems 2.1.7 and 2.1.9]{Mu}),
the Gelfand spectrum of $A$ is identified with the set of *-homomorphic functionals which are $\| \cdot\|_1$-bounded, i.e.,
$(\varphi(h_n))$ being bounded as an element in $(\ell^1)^* = \ell^\infty$.  
In view of the expression
  \[
    \varphi(h_n) = \alpha \gamma_+^n + \beta \gamma_-^n,
    \quad
 \gamma_\pm = \frac{c \pm \sqrt{c^2 - 4r(1-r)}}{2(1-r)},  
  \]
$\sup_{n \geq 0} |\varphi(\delta_n)| < \infty$ if and only if $|\gamma_\pm| \leq 1$.
The last condition is then rephrased, by a simple argument, in terms of $c \in \R$ as $-1 \leq c \leq 1$,
which remains valid even for degenerate cases $c = \pm 2\sqrt{r(1-r)}$.
% If $c = \pm 2\sqrt{r(1-r)}$, $\phi(\delta_n) = \left( \pm \sqrt{\frac{r}{1-r}} \right)^n( 1+ (1-2r)n)$ ($n \geq 1$),
% which is bounded for $0 \leq r \leq 1/2$. 

As a summary, the spectrum of $h_1$ in $A$ is the interval $[-1,1]$.
% The range assumption $0 \leq r \leq 1/2$ is necessary here. 

% As a polynomial algebra, a multiplicative functional $\varphi$
% on $\sD$ is given as an evaluation of $t$ at some (and any) $\gamma \in \C$;
% $\varphi$ is specified by the value $\gamma = \varphi(\delta_1)$ and
% $\varphi(\delta_n)$ is determined recursively by the relation
% $\gamma \varphi(\delta_n) = r\varphi(\delta_{n-1}) + (1-r) \varphi(\delta_{n+1})$ ($n \geq 1$) 
% together with the initial value $\varphi(\delta_0) = 1$.

{\small
\begin{Remark}
  Similarly bounded homomorphisms are characterized by an ellipse condition on $c \in \C$; $c = a + ib$ ($a,b \in \R$) satisfying
  $a^2 + (b/(1-2r))^2 \leq 1$ 
(see \cite{Ca,Py,PS} and \cite[\S 3.3]{FP} for this and related results). 
\end{Remark}}

  \section{Analytic Functionals}
  Consider a linear functional $\varphi$ on $\sA$ which is \textbf{analytic} in the sense that the power series
  $\varphi(z) = \sum_{n \geq 0} \varphi(h_n) z^n$ is convergent in a neighborhood of $z=0$.
  We shall investigate spectral properties of $\varphi$ by computing its Stieltjes transform. 
 % inversion formula of the Cauchy transform of $h$. 
  
% Given a real $s>0$, we introduce a hermitian functional $h$ (called Haagerup functional after \cite{Ha})
% on $\sA$ by $h(\delta_n) = e^{-sn}$ ($n \geq 0$).
% We shall show $h$ is positive based by solving the accompanied moment problems in an explicit way.

Consider a generating function
$p(z,t) = \sum_{n=0}^\infty z^n p_n(t)$ of the polynomial sequence $( p_n(t) )_{n \geq 0}$.
From the algebraic relation 
\begin{align*}
  tp(z,t) &= t + \sum_{n=1}^\infty rz^n p_{n-1}(t) + (1-r) \sum_{n=1}^\infty z^n p_{n+1}(t)\\
          &= t + rz p(z,t) + \frac{1-r}{z} (p(z,t) - 1 - zt), 
\end{align*}
we obtain a closed expression 
\[
  p(z,t) = \frac{1-r-rzt}{1-r -zt + rz^2}.
\]
The evaluation by $\varphi$ (which is formally interpreted as an integration with respect to a symbolical measure
$\mu(dt)$ on $\R$) gives
\[
  \int \frac{1-r-rzt}{1-r -zt + rz^2}\, \mu(dt)
  = \sum_{n=0}^\infty z^n \varphi(h_n)
  = \varphi(z). 
%  = \frac{1}{1- e^{-s}z}. 
\]
By rewriting the integrand, we have
\[
%  \frac{1}{1- e^{-s}z}
  \varphi(z) = r\varphi_0 + \frac{r^2z^2 - (1-r)^2}{z} \int \frac{1}{t - rz - \frac{1-r}{z}}\, \mu(dt). 
\]
The case $r=0$ is simple and separately dealt with: In terms of a parameter $w = 1/z$, this takes the form 
\[
\int \frac{1}{t-w}\, \mu(dt) = - \frac{1}{w} \varphi\bigl(\frac{1}{w}\bigr), %\frac{1}{e^{-s} - w}
\]
which is combined with Theorem~\ref{characterization} to see that the functional $\varphi$ belongs to $A^*$
(i.e., $\mu$ is a genuine complex Radon measure on $\R$) if and only if
$\varphi(1/w)$ is holomorphically extended to $\C \setminus [-1,1]$ in such a way that
% \[
%   \lim_{y \to \pm 0} \frac{1}{t+iy} \varphi\left(\frac{1}{t+iy}\right)
% \]
% exist as distributions on $C_c^\infty(\R)$ and 
the Stieltjes inversion formula
\[
  2\pi i\mu(dt) = \lim_{\epsilon \to +0} \left( \frac{1}{t-i\epsilon} \varphi \Bigl(\frac{1}{t-i\epsilon}\Bigr)
  - \frac{1}{t+i\epsilon} \varphi \Bigl(\frac{1}{t+i\epsilon}\Bigr) \right) dt
\]
gives a complex measure $\mu$ on $[-1,1]$. 
% and the functional $h$ is given by a Dirac measure located at $t = e^{-s}$.

From here on we assume $0 < r \leq 1/2$ and introduce a new variable $w$ by 
\[
  w = rz + \frac{1-r}{z} \iff  z = \frac{w - \sqrt{w^2 - 4r(1-r)}}{2r} 
\]
so that the Stieltjes transform $S(w)$ of $\varphi$ (symbolically $\mu$) is given by
\[
  \int \frac{1}{t-w}\, \mu(dt) = (\varphi(z) - r\varphi_0) \frac{z}{r^2z^2 - (1-r)^2}.
%  = \frac{1-r + e^{-s}rz}{r - e^{-s}rz} \frac{rz}{r^2z^2 - (1-r)^2}.  
\]

Here the square root is chosen so that the right hand side behaves
like $-\varphi_0/w - \varphi_1/w^2 + \cdots$, i.e., $\sqrt{w^2 - 4r(1-r)} \approx w$ or equivalently
$z \approx (1-r)/w$ for a large $w$.
Note that as a holomorphic function of $w$, $\sqrt{w^2 - 4r(1-r)}$ is defined
on $\C \setminus \sigma_r$ with
\[
  \sigma_r = [-2\sqrt{r(1-r)},2\sqrt{r(1-r)}] \subset [-1,1]
\]
and continuously extended to its closure
$\C \setminus (-2\sqrt{r(1-r)},2\sqrt{r(1-r)})$.
Even on the cut line, it approaches continuous functions separately:
\[
  \lim_{s \to \pm 0} \sqrt{(t+is)^2 - 4r(1-r)} =
 \begin{cases} \pm i \sqrt{4r(1-r) - t^2} &
   (t^2 \leq 4r(1-r)),\\
   \sqrt{t^2 - 4r(1-r)} & (t \geq 2\sqrt{r(1-r)},\\
   - \sqrt{t^2 - 4r(1-r)} & (t \leq -2\sqrt{r(1-r)}.
 \end{cases}
\]

% In the limit $s \to \infty$,
When $\varphi(z) \equiv 1$, the above formula is reduced to
\[
  \int \frac{1}{t-w}\, \mu(dt) = \frac{(1-r)z}{r^2z^2 - (1-r)^2}
  = \frac{-(1-2r)w + \sqrt{w^2 - 4r(1-r)}}{2r(1-w^2)} 
\]
and an explicit computation of the Stieltjes inversion formula reveals that, for $0 < r \leq 1/2$, $\mu$ is 
a continuous measure 
\[
  \frac{1}{2\pi r} \frac{\sqrt{4r(1-r) - t^2}}{1-t^2} dt
\]
supported by the interval $\sigma_r$ with no atomic measures arising from $w= \pm 1$. 
% , for $1/2 < r < 1$,adding to this an atomic measure
% \[
% \frac{2r-1}{2r} (\delta(t-1) + \delta(t+1)) 
% \]
% of negative weights appears.

Thus $\varphi(z) \equiv 1$ gives rise to a state (a probability measure supported by $[-1,1]$) of $A$ as a linear functional. 
% Thus the linear functional $h$ in the limit $s \to \infty$ is positive.
% if and only if $0 < r \leq 1/2$.
Notice that at the boundary parameter $r = 1/2$ this is further reduced to 
\[
  \mu(dt) = \frac{1}{\pi} \frac{1}{\sqrt{1-t^2}}\, dt.
\]

{\small
\begin{Remark}
  The measure for $\varphi(z) \equiv 1$ first appeared in \cite{Ke} and 
  was used as a Plancherel measure in \cite{Ca} by a similar method presented here 
  (see \cite[\S 3.4]{FP} for details). 
\end{Remark}}

Returning to the general case, introduce a function $S(w)$ of $w$ by 
\[
S(w) =   \frac{(2r-1)w + \sqrt{w^2 - 4r(1-r)}}{2r(1-r)(1-w^2)}\ 
   \left(
     \varphi \Bigl(
      \frac{w - \sqrt{w^2 - 4r(1-r)}}{2r} 
      \Bigr) - r\varphi_0
      \right), 
\]
which is first defined for a large $w$ and expanded at $w = \infty$ as 
\[
  S(w) = -\frac{\varphi_0}{w}  - \frac{\varphi_1}{w^2} - \frac{r\varphi_0 + (1-r)\varphi_2}{w^3}
  - \frac{r(2-r)\varphi_1 + (1-r)^2 \varphi_3}{w^4} + \cdots.
\]

These are now combined with Theorem~\ref{characterization} to have the following. 

\begin{Theorem}
  An analytic functional $\varphi$ on $\sA$ belongs to $A^*$ if and only if
  \begin{enumerate}
  \item the function $S$ is holomorphically extended to $\C \setminus [-1,1]$ and 
  % \item
  %   limits $\displaystyle \lim_{y \to \pm 0} S(x+iy)$ 
  %   exist as distributions on $C_c^\infty(\R)$ and
\item
  the Stieltjes inversion formula
\[
  2\pi i\mu(dt) = \lim_{\epsilon \to +0} \left(
 S(t+i\epsilon) - S(t-i\epsilon) 
  \right) dt
\]
gives a complex Radon measure $\mu$ on $[-1,1]$.
\end{enumerate}
\end{Theorem}

{\small
\begin{Remark}
  By reversing the process, we see that any functional which belongs to $A^*$ is analytic: 
  As a Stieltjes transform of a complex measure supported by $[-1,1]$, $S(w)$ is analytic at $w=\infty$,
  which in turn implies the analyticity of $\varphi(z)$ at $z = 0$. 
\end{Remark}}

\section{Analytic Functionals of Geometric Series}
Let $\varphi(z) = \varphi_\upsilon(z) = 1/(1-\upsilon z)$ with $\upsilon \in \C$, i.e.,
$\varphi_n = \varphi(h_n) = \upsilon^{n}$ ($n \geq 1$) and $\varphi_0 = 1$.
Since the case $\upsilon = 0$, i.e., $\varphi(z) \equiv 1$ is already discussed, we assume $\upsilon \not= 0$ in what follows. 

First consider the case $r=0$. Then $S(w) = 1/(\upsilon - w)$, which is holomorphic on $\C \setminus [-1,1]$ if and only if
$-1 \leq \upsilon \leq 1$. In that case, $S$ is the Stieltjes transform of a Dirac measure at $t = \upsilon \in [-1,1]$.

When $0 < r \leq 1/2$, 
\[
  2(\varphi(z) - r\varphi_0) % = 2 \frac{1-r + r\upsilon z}{1-\upsilon z}
  = \frac{2r(1-r)(\upsilon^{-1} - \upsilon) + (2r-1)w - \sqrt{w^2 - 4r(1-r)}}{c_r(\upsilon) - w}.  
\]
and we have 
% In the numerator of $2(h(z) - rh_0)$, singular points are $w = \pm 2\sqrt{r(1-r)}$,
% one of which coincides with $c_r(\upsilon)$ if and only if $\upsilon^2 = (1-r)/r$.
% Thus, for $\upsilon \not= \pm \sqrt{(1-r)/r}$, we can evaluate the numerator in $2(h(z) - rh_0)$
% by putting $w = c_r(\upsilon)$ to be
% \[
%   r\upsilon - (1-r)\upsilon^{-1} - \sqrt{c_r(\upsilon)^2 - 4r(1-r)}. 
% \]
% % with $c_r(\upsilon)^2 - 4r(1-r) = (r\upsilon - (1-r)\upsilon^{-1})^2$.
% Since we have
% \[
%   \sqrt{c_r(\upsilon)^2 - 4r(1-r)} = r\upsilon - (1-r)\upsilon^{-1}
% \]
% for a large positive real
% $\upsilon$ and the numerator at $w = c_r(\upsilon)$
% is locally holomorphic in $\upsilon \in \C \setminus\{ 0, \pm \sqrt{(1-r)/r}\}$, it vanishes identically.
% In this way, we have seen that
% $h(z) - rh_0$ is holomorphic as a function of $w$ in a neighborhood of $c_r(\upsilon)$
% if $\upsilon^2 \not= (1-r)/r$.
% Now we put the things together and have
\[
  S(w) = \frac{1}{c_r(\upsilon) - w}
  + \frac{\upsilon^{-1} - \upsilon}{2} \frac{(2r-1)w + \sqrt{w^2 - 4r(1-r)}}{(1-w^2)(c_r(\upsilon) - w)}, 
\]
where $c_r(\upsilon) = r\upsilon^{-1} + (1-r)\upsilon$ is a weighted Joukowsky transform of $\upsilon$ %$c_r((1-r)/r\upsilon)$ 
and gives a possible point of singularity other than $w = \pm 1$ and $w = \pm 2\sqrt{r(1-r)}$.
With a polar expression $\upsilon = |\upsilon| e^{i\theta}$, if $\theta$ moves from $0$ to $2\pi$, 
\[
  c_r(\upsilon) = (r|\upsilon|^{-1} + (1-r)|\upsilon|)\cos\theta + i(-r|\upsilon|^{-1} + (1-r)|\upsilon|) \sin\theta
\]
traces an ellipse surrounding the interval $\sigma_r$ 
clockwise or counterclockwise according to $|\upsilon| > \sqrt{r/(1-r)}$ or
$|\upsilon| < \sqrt{r/(1-r)}$. At the critical value $|\upsilon| = \sqrt{r/(1-r)}$, 
the ellipse shrinks to the interval $\sigma_r$. 

{\small
\begin{Remark}
  As a real-valued function of $\upsilon \geq 0$, $c_r(\upsilon)$ attains the minimal value $2\sqrt{r(1-r)}$
  at $\upsilon = \sqrt{r/(1-r)}$ and takes a value $1$
exactly at $\upsilon = r/(1-r)$ or at $\upsilon = 1$.
\end{Remark}}

From 
$c_r(\upsilon)^2 - 4r(1-r) = (r\upsilon^{-1} - (1-r)\upsilon)^2$,
one sees that 
\[
  \sqrt{w^2 - 4r(1-r)}|_{w=c_r(\upsilon)} = \pm (r\upsilon^{-1} - (1-r)\upsilon)
\]
with $\pm$ determined as follows:
When $\theta = 0$, $c_r(\upsilon) = c_r(|\upsilon|) \geq 2\sqrt{r(1-r)}$ and we have 
$\sqrt{c_r(|\upsilon|)^2 - 4r(1-r)} = \bigl|r|\upsilon|^{-1} - (1-r)|\upsilon|\bigr|$, i.e., 
$r|\upsilon|^{-1} - (1-r)|\upsilon|$ or $(1-r)|\upsilon| - r|\upsilon|^{-1}$
according to $|\upsilon| < \sqrt{r/(1-r)}$ or $|\upsilon| > \sqrt{r/(1-r)}$,
which is then analytically continued to the identity 
\[
  \sqrt{c_r(\upsilon)^2 - 4r(1-r)} 
=
\begin{cases}
  r\upsilon^{-1} - (1-r)\upsilon &\text{if $|\upsilon| < \sqrt{r/(1-r)}$,}\\
  (1-r)\upsilon - r\upsilon^{-1} &\text{if $|\upsilon| > \sqrt{r/(1-r)}$.}
\end{cases}
\]
Consequently we have 
\begin{multline*}
(2r-1) c_r(\upsilon)+ \sqrt{c_r(\upsilon)^2 - 4r(1-r)}\\ 
=
\begin{cases}
  2r^2\upsilon^{-1} - 2(1-r)^2\upsilon &\text{if $|\upsilon| < \sqrt{r/(1-r)}$,}\\
  2r(1-r)(\upsilon - \upsilon^{-1}) &\text{if $|\upsilon| > \sqrt{r/(1-r)}$}
\end{cases}
\end{multline*}
and, after a simple computation, 
\[
  \lim_{w \to c_r(\upsilon)} (c_r(\upsilon) - w) S(w) =
   \begin{cases}
    0 &\text{if $|\upsilon| < \sqrt{r/(1-r)}$,}\\
    \frac{1-c_r(\upsilon^2)}{1-c_r(\upsilon)^2} &\text{if $|\upsilon| > \sqrt{r/(1-r)}$.}
  \end{cases}
\]
% \[
%   \lim_{w \to c_r(\upsilon)} (c_r(\upsilon) - w) C(w) =
%   \frac{2(1-c^2) + (\upsilon - \upsilon^{-1})((2r-1)c \pm (r\upsilon - (1-r)\upsilon^{-1}))}{2(1-c^2)}
% \]
% \[
%   C(w) = \frac{2(1-w^2) + (2r-1)(\upsilon - \upsilon^{-1})w + (\upsilon - \upsilon^{-1})
%     \sqrt{w^2 - 4r(1-r)}}{2(1-w^2)(c_r(\upsilon) - w)}, 
% \]
% where $w = c_r(\upsilon)$ is a fake singularity unless $\upsilon^2 = (1-r)/r \iff c_r(\upsilon)^2 = 4r(1-r)$.

Thus an atomic measure arises at the point $c_r(\upsilon) \not\in \sigma_r$ only if $|\upsilon| > \sqrt{r/(1-r)}$
with the weight given by $(1-c_r(\upsilon^2))/(1-c_r(\upsilon)^2)$.
Note that $c_r(\upsilon) \in [-1,1] \setminus \sigma_r^\circ$
if and only if $r/(1-r) \leq |\upsilon| \leq 1$ with $\upsilon \in \R$.

There remains the possibility of an atomic measure at points $w=\pm 1$ arising from the second part of $S$ but this is not the case
in view of 
\[
  \lim_{w \to \pm 1} ((2r-1)w + \sqrt{w^2 - 4r(1-r)})
  = \pm(2r-1) \pm (1-2r) = 0. 
\]

We next look into the measure supported by $\sigma_r$. % which is necessarily supported by the interval $I_r$.

(i) $c_r(\upsilon) \in \sigma_r^\circ \iff |\upsilon| = \sqrt{r/(1-r)}$ and $\upsilon \not\in \R$:
$S(t+is) - S(t-is)$ exhibits a singularity proportional to $1/(c_r(\upsilon) - t)$ inside $\sigma_r$ and
the Stieltjes inversion formula fails to give a measure.

(ii) $c_r(\upsilon) = \pm 2\sqrt{r(1-r)} \iff \upsilon = \epsilon \sqrt{r/(1-r)}$ with $\epsilon = \pm 1$:
  From the relation
  \[
    \lim_{w \to c_r(\upsilon)} \frac{\upsilon^{-1} - \upsilon}{2} \frac{(2r-1)w}{1-w^2} = -1, 
  \]
  the inversion formula for the part
  \[
    \frac{1}{c_r(\upsilon) - w} + \frac{\upsilon^{-1} - \upsilon}{2} \frac{(2r-1)w}{(1-w^2)(c_r(\upsilon) - w)}
  \]
  gives zero contribution as a measure.

  For the continuous part, we observe that
  \[
  \frac{\sqrt{|(t+ is)^2 - 4r(1-r)|}}{|c_r(\upsilon) - (t+ is)|} =  
  \sqrt{\frac{|c_r(\upsilon) + t + is|}{|c_r(v) - (t+is)|}} 
  \leq \sqrt{\frac{|c_r(v) + t| + |s|}{|c_r(\upsilon) - t| \vee |s|}}
\]
and % the maximum is in fact attained. 
\[
  \sup_{t \in \R} \frac{|c_r(v) + t| + |s|}{|c_r(\upsilon) - t| \vee |s|}
  = \frac{2(|s| + |c_r(\upsilon)|)}{|s|} \vee \frac{|s|}{2|c_r(\upsilon)|}, 
\]
whence
\[
  \frac{\sqrt{|(t+ is)^2 - 4r(1-r)|}}{|c_r(\upsilon) - (t+ is)|} \leq \sqrt{\frac{|c_r(\upsilon) + t| + |s|}{|c_r(\upsilon) - t|}}, 
\]
are dominated by a locally integrable function $\sqrt{\frac{|c_r(\upsilon) + t| + 1}{|c_r(\upsilon) - t|}}$ of $t \in \R$
uniformly in $0 < |s| \leq 1$ and 
\[
  \lim_{s \to 0} |s| \frac{\sqrt{|(t+ is)^2 - 4r(1-r)|}}{|c_r(\upsilon) - (t+ is)|} = 0
 % \leq |s| \sqrt{\frac{|c_r(v) + t| + |s|}{|c_r(\upsilon) - t| \vee |s|}},  
\]
uniformly in $t \in \R$. 

Now, in view of the expression 
  \[
    \frac{\sqrt{(t+ is)^2 - 4r(1-r)}}{c_r(\upsilon) - (t+ is)}
    = i \frac{s}{|s|} \frac{\sqrt{2\sqrt{r(1-r)} + \epsilon (t+is)}}{\sqrt{2\sqrt{r(1-r)} - \epsilon (t+is)}}
       \quad
    (t \in \sigma_r), 
  \]
 the dominated convergence theorem is applied to have 
  % \[
  %   \frac{(16r(1-r) + s^2)^{1/4}}{|s|^{1/2}}
  % \]
  % for $0 \not= s \in \R$, whence we see
\begin{multline*}
  \lim_{s \to \pm 0} (\upsilon^{-1} - \upsilon)
  \frac{\sqrt{(t+ is)^2 - 4r(1-r)}}{c_r(\upsilon) - (t+ is)}\,dt\\
  = \pm i \left( \sqrt{\frac{1-r}{r}} - \sqrt{\frac{r}{1-r}} \right)
  \sqrt{\frac{2\sqrt{r(1-r)} + \epsilon t}{2\sqrt{r(1-r)} - \epsilon t}}\,dt 
\end{multline*}
as a complex Radon measure on $-2\sqrt{r(1-r)} < t < 2\sqrt{r(1-r)}$. 

Thus a positive finite measure supported by $\sigma_r$ comes out as their difference divided by
$2\pi i$. 
% Notice here that, in view of
% \[
%   \frac{\sqrt{|(t+ is)^2 - 4r(1-r)|}}{|c_r(\upsilon) - (t+ is)|} =  
%   \frac{\sqrt{|c_r(v) + (t+is)|}}{\sqrt{|c_r(v) - (t+is)|}}
%   \leq \frac{\sqrt{|c_r(v) + (t+is)|}}{\sqrt{|2\sqrt{r(1-r)} - \epsilon t|}}, 
% \]
% the weight functions of the measure for $0 < |s| \leq 1$ are dominated by an integrable function of $t \in [-1,1]$. 

% Thus both limits $s \to \pm 0$ exist as locally finite continuous measures thanks to the dominated convergence theorem 
% and a continuous measure supported by $\sigma_r$ comes out as their difference. 
% supported by $(-2\sqrt{r(1-r)},2\sqrt{r(1-r)})$ appears.
% is, including the case $\upsilon^2 = (1-r)/r$,

(iii) $c_r(\upsilon) \not\in \sigma_r \iff |\upsilon| \not= \sqrt{r/(1-r)}$:
In view of 
\[
    \lim_{s \to \pm 0} (\upsilon^{-1} - \upsilon)
    \frac{\sqrt{(t+ is)^2 - 4r(1-r)}}{c_r(\upsilon) - (t+ is)}
    = \pm i (\upsilon^{-1} - \upsilon)  \frac{\sqrt{4r(1-r) - t^2}}{c_r(\upsilon) - t} 
  \]
for $t \in \sigma_r$, a continuous measure supported by $\sigma_r$ appears again and its explicit form is given by the theorem below, 
% \[
%   \frac{\upsilon - \upsilon^{-1}}{2\pi}
%   \frac{\sqrt{4r(1-r) - t^2}}{(1-t^2)(c_r(\upsilon) - t)}\, dt,  
% \]
which covers the measure in (ii) if one puts $\upsilon = \pm \sqrt{r/(1-r)}$.

{\small
\begin{Remark}
The continuous measure in (ii) or (iii)
is real if and only if both $\upsilon^{-1} - \upsilon$ and $c_r(\upsilon)$ are real, i.e.,
if $\upsilon \in \R$. 
Furthermore it is positive if and only if $\upsilon^{-1} - \upsilon$ and $c_r(\upsilon)$ have same signatures,
i.e., $-1 \leq \upsilon \leq 1$. % $\upsilon \leq -1$ or $\upsilon \geq 1$. 
\end{Remark}}

\begin{figure}[h]%{r}[overhang]{width}
  \centering
 \includegraphics[width=0.7\textwidth]{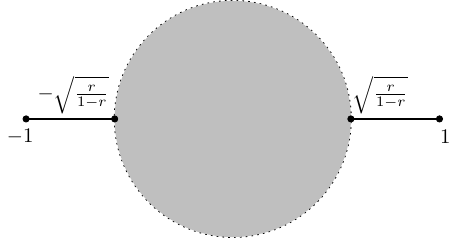}
 \caption{$\upsilon$-region of C*-boundedness}
\end{figure}

Summarizing discussions so far, we have the following.

\begin{Theorem}\label{main}
 Let $0 < r \leq 1/2$. The analytic linear functional of geometric series parametrized by $\upsilon \in \C$ belongs to $A^*$ if and only if
  the parameter $\upsilon$ falls into one of the following cases. 
  % $1 \leq |\upsilon|$ with $\upsilon \in \R$.
  \begin{enumerate}
    \item
  For $|\upsilon| < \sqrt{r/(1-r)}$ or $\upsilon = \pm \sqrt{r/(1-r)}$, the associated measure is 
  continuous, supported by the interval $\sigma_r$ and of the form 
  \[
     \frac{\upsilon^{-1} - \upsilon}{2\pi}
  \frac{\sqrt{4r(1-r) - t^2}}{(1-t^2)(c_r(\upsilon) - t)}\, dt.   
\]
\item 
For a real $\upsilon$ satisfying $\sqrt{r/(1-r)} < |\upsilon| \leq 1$, adding to this continuous measure,
there appears an atomic measure of the form
\[
  \frac{1-c_r(\upsilon^2)}{1-c_r(\upsilon)^2} \delta(t-c_r(\upsilon))
  = \frac{(1-r)\upsilon^2 - r}{(1-r)^2\upsilon^2 - r^2} \delta(t-c_r(\upsilon)). 
\]
\end{enumerate}

Among these, positive ones are given by $-1 \leq \upsilon \leq 1$. % $\pm \upsilon \leq 1$. 
\end{Theorem}

{\small
\begin{Remark}~
  \begin{enumerate}
  \item
    Under the transform $\upsilon \mapsto c_r(\upsilon)$,
    an open disk $|\upsilon| < \sqrt{r/(1-r)}$ is biholomorphically
    mapped to a slitted Riemann sphere $\overline{\C} \setminus \sigma_r$. 
  \item
    Both $c_r(\upsilon)^2$ and $c_r(\upsilon^2)$ are strictly positive and less than $1$ with
    $c_r(\upsilon)^2 = 1 \iff \upsilon = \pm 1$ and $c_r(\upsilon^2) = 1 \iff \upsilon^2 \in \{1,r/(1-r)\}$.
   When $\upsilon= \pm 1$, the coefficient in (ii) of Theorem takes an indefinite form and 
%Note that $c_r(\upsilon) = \pm 1$ for $1 \leq |\upsilon| < (1-r)/r$ if and only if $\upsilon = \pm 1$.
    % This exceptional situation gives the value $1$ and
    should be dealt with directly: $S(w) = 1/(c_r(\upsilon) -w)$ shows that it gives a Dirac measure
    at $t= \pm 1$ and the formula can be understood $(1-c_r(\upsilon^2)/(1-c_r(\upsilon)^2) = 1$
    in this case.
  \item
  When a positive $\upsilon$ increases and passes the point $\sqrt{r/(1-r)}$,
    the same form of continuous measure reappears with a constant less than $1$ multiplied and
    an atomic measure appears to compensate the deficient mass. At the critical point $v = \sqrt{r/(1-r)}$ there is still no atomic measure
    but the density function of the continuous measure diverges at a spectral boundary $2\sqrt{r(1-r)}$ to release a point mass
    at $c_r(v)$ for $v > \sqrt{r/(1-r)}$. (Notice that the density function vanishes at $2\sqrt{r(1-r)}$ for $v \not= \sqrt{r/(1-r)}$.)
    After the critical point, the continuous mass decreases and vanishes at $v=1$, leaving an atomic probability at a spectral end $t=1$. 
  \item
    For a real $\upsilon$ satisfying $\sqrt{r/(1-r)} < |\upsilon| \leq 1$, let $\upsilon' = \frac{r}{1-r} \frac{1}{\upsilon}$ so that
    $c_r(\upsilon) = c_r(\upsilon')$ and $\frac{r}{1-r} \leq |\upsilon'| < \sqrt{r/(1-r)}$. Then
    \[
      \frac{1-c_r(\upsilon^2)}{1- c_r(\upsilon)^2} + \frac{\upsilon^{-1} - \upsilon}{(\upsilon')^{-1} - \upsilon'} = 1. 
    \]
    and the measure for $\upsilon$ is a convex combination of the Dirac measure at $c_r(\upsilon)$ and the measure for $\upsilon'$. 
  % \item
  % By a partial fraction decomposition
  % $P(z,t) = b_\upsilon \frac{\upsilon}{\upsilon - z} + b_{\upsilon'} \frac{\upsilon'}{\upsilon' - z}$ with
  % $\upsilon\upsilon' = (1-r)/r$ and $\upsilon + \upsilon' = t/r$, one sees that the evaluation functional
  % at $c_r(\upsilon)$ is equal to $b_\upsilon h_\upsilon + b_{\upsilon'} h_{\upsilon'}$, where
  % \[
  %   b_\upsilon = \frac{\upsilon^2 - (1-r)^2/r^2}{\upsilon^2 - (1-r)/r}.
  % \]
      \item
  We explicitly computed the Stieltjes inversion formula to get a finite bounded measure,
  which in turn gives rise
  to the starting holomorphic function as a Stieltjes transform. In fact, both $S(w)$ and
  $\int \frac{1}{t-w}\, \mu(dt)$ have the same boundary value on the real line, whence their difference
  is entirely holomorphic and turns out to be identically zero in view of their behavior at infinity.
  (See Appendix~B for more information.)
% \item
%  Taking $\upsilon \to 0$, the measure converges to the Kesten measure.
% \item
%   When $r = 1/2$, the theorem is reduce to the following:
%  The analytic linear functional of geometric series parametrized by $\upsilon \in \C$ belongs to $A^*$ if and only if
%   the parameter $\upsilon$ falls into one of the following two cases. 
%   % $1 \leq |\upsilon|$ with $\upsilon \in \R$.
%   \begin{enumerate}
%     \item
%   For $|\upsilon| < 1$, the associated measure is 
%   continuous, supported by the interval $[-1,1]$ and of the form 
%   \[
%    \frac{1}{2\pi}  \frac{\upsilon^{-1} - \upsilon}{((c_{1/2}(\upsilon) - t)\sqrt{1-t^2}}\, dt,   
%  \]
%  which is reduced to the arcsine law for $\upsilon = 0$. 
% \item
%   For $\upsilon = \pm 1$, the associated measure is $\delta(t\mp 1)$, which is the limit of (a) as $\upsilon \to \pm 1$. 
% \end{enumerate}
% 
% Among these, positive ones are given by $-1 \leq \upsilon \leq 1$. % $\pm \upsilon \leq 1$. 
\end{enumerate}
\end{Remark}}

\begin{Corollary}[Normalized Limit]
   In terms of the normalized variable $\displaystyle s = \frac{t-\upsilon}{\sqrt{r(1-\upsilon^2)}}$,
    the normalized distribution converges to $\nu(ds)$
      under the limit $(r,\upsilon) \to (0,0)$ keeping $\gamma = \upsilon/\sqrt{r} \in \R$ constant (see \cite[\S 4.3]{HO} for more information):
      \begin{enumerate}
      \item
        For $-1 \leq \gamma \leq 1$,
        \[
          \nu(ds) = \frac{ \sqrt{4 - (s+\gamma)^2}}{2\pi(1-\gamma s)}\, ds.
        \]
      \item
        For $\pm\gamma \geq 1$,
        \[
          \nu(ds) =  \frac{ \sqrt{4 - (s+\gamma)^2}}{2\pi(1-\gamma s)}\, ds
          + \left(1-\frac{1}{\gamma^2}\right) \delta(s - 1/\gamma). 
        \]
      \end{enumerate}
\end{Corollary}

\begin{Corollary}
  In terms of a new variable $x = (1-\gamma s)/\gamma^2$, the normalized limit becomes the free Poisson distribution of mean 
  $\lambda = 1/\gamma^2$ (see \cite{MS} for example as a reference):
  \begin{enumerate}
  \item
    For $\lambda \geq 1$,
    \[
     \frac{\sqrt{4\lambda - (x-1-\lambda)^2}}{2\pi x}\, dx. 
      % \quad
      % ((1-\lambda)^2 \leq x \leq (1+\lambda)^2). .
    \]
  \item
    For $0 < \lambda < 1$,
    \[
      \frac{\sqrt{4\lambda - (x-1-\lambda)^2}}{2\pi x}\, dx +
      (1 - \lambda) \delta(x).
    \]
  \end{enumerate}
\end{Corollary}

{\small
\begin{Remark}
  The free Meixner distribution (see \cite{BB,SY}) also contains the free Poisson distribution as a special case
  and it is an interesting problem to merge these with the distribution in Theorem 3.1.
  For this, \cite{YY} might be helpful. 
\end{Remark}
}

\section{Positive Definite  Functions of Haagerup}
% \begin{Remark}
%   The case $ r = 1/2l$ with $l \geq 1$ a natural number is associated to a free group of $l$ generators. 
% \end{Remark}

We shall here explain how the positive functionals considered in the previous section are related to
positive definite functions of Haagerup on free groups (\cite{Ha,Bo}). 
% As an application in the previous section, we give a proof of the fact that  Haagerup functionals
% positive definite on free groups. The original proof due to Haagerup depends heavily on
% combinatoral analysis in free groups.

Let $G$ be a free group of $l$ generators and $G_n$ be the set of words of length $n$ so that
$G = \sqcup_{n \geq 0} G_n$. Note that $|G_n| = 2l(2l-1)^{n-1}$ ($n \geq 1$). 
A radial function is by definition a function of $g \in G$ whose values depend only on the word length $|g|$.
Let $\sA$ be the set of elements in the group algebra $\C G$ consisting of radial functions.
As a linear basis of $\sA$, we choose $h_n = \frac{1}{|G_n|} 1_{G_n}$ ($n \geq 0$).
It then satisfies $h_1 h_n = \frac{1}{2l} h_{n-1} + \frac{2l-1}{2l} h_{n+1}$ and
$\sA$ is a commutative *-subalgebra of $\C G$.

Let the averaging map $E: \C G \to \sA$ be defined by $E(g) = h_{|g|}$ ($g \in G \subset \C G$), 
which satisfies the algebraic properties of conditional expectation:
$E(a) = a$, $E(f^*) = E(f)^*$ and $E(af) = aE(f)$ for $a \in \sA$ and $f \in \C G$.

On the regular representation of $\C G$,
$E$ is `implemented' by the projection $e$ to the closed subspace $\overline{\sA \tau^{1/2}}$
($\tau$ being the standard trace on $C^*(G)$ with $\tau^{1/2}$ the GNS vector) in the sense that
$E(f)e = e f e$ ($f \in \C G$) and hence
induces by continuity a conditional expectation (also denoted by $E$) on the group von Neumann algebra in such a way that
$E(C^*_{\text{red}}(G)) \subset C^*_{\text{red}}(G)$. 
Here $C^*_\text{red}(G)$ is the reduced group C*-algebra, i.e., the C*-algebra generated by the regular representation of $\C G$.

According to Haagerup, we introduce a positive definite function of $G$ by 
$\phi_\upsilon(g) = \upsilon^{|g|}$ ($g \in G \subset \C G$) for $0 < \upsilon < 1$ with the associated state 
on the group C*-algebra $C^*(G)$ also denoted by $\phi_\upsilon$.
Notice that $\phi_\upsilon(f) = \varphi_\upsilon(E(f))$ for $f \in \C G$. 
Thus, when the state $\phi_\upsilon$ % of $C^*(G)$ associated with a Haagerup's positive definite function $h$ 
splits through $C^*_\text{red}(G)$,
it is induced from the restricted state on a commutative C*-subalgebra $A_{\text{red}} = E(C^*_\text{red}(G))$ of $C^*_\text{red}(G)$
via the conditional expectation $C^*_\text{red}(G) \to A_\text{red}$.
This is exactly the case when $\upsilon \leq \sqrt{r/(1-r)}$ ($r = 1/2l$) as seen below.

Since $\tau$ is a limit case $\displaystyle \lim_{\upsilon \to 0}\phi_\upsilon$,
the analysis in \S 2 shows that the spectrum of $h_1$ in $C^*_\text{red}(G)$ is
the interval $\sigma_r$ as already seen in \cite{Ke}.
% Also we can show that any C*-seminorm on $\sA$ is continuous with respect to the C*-norm of $C^*(G)$. 

In view of the result in the previous section, 
$\phi_\upsilon$ is expected to be positive definite even for $-1 < \upsilon < 0$. This is in fact the case because
$g \mapsto (-1)^{|g|} g$ ($g \in G$) induces an involutive *-automorphism  of $\C G$. 
As a consequence, we observe that $\sA \to C^*(G)$ is extended to an embedding $A \to C^*(G)$ of a C*-algebra.
Non-trivial is injectivity,
which holds because we have a family of states $(\phi_\upsilon) _{-1 \leq \upsilon \leq 1}$ of $C^*(G)$ with
their pullbacks to $A$ providing a faithful family. 
%We know that $\| a\| = \sup\{ |h_\upsilon(a)|; -1 \leq \upsilon \leq 1\}$ for $a \in \sA$.
%$\sup\{ |h_\upsilon(a)|; \upsilon \geq \sqrt{r(1-r)}\} = \sup\{ |a(t)|; t \in I_r\}$ for $a \in \sA$.
Consequently the spectrum of $h_1$ in $C^*(G)$ is equal to that in $A$, i.e., $[-1,1]$,
and we shall identify the radial subalgebra $A$ of $C^*(G)$ with $C([-1,1])$ in the remaining. % as a unital C*-subalgebra of $C^*(G)$.

From these observations, the linear functional $\phi_\upsilon$ is continuous on $C^*_\text{red}(G)$ if and only if
the spectral measure of $\varphi_\upsilon = \phi_\upsilon|_A$ is supported by $\sigma_r$, i.e.,
$|\upsilon| \leq \sqrt{r/(1-r)} = 1/\sqrt{2l-1}$, which agrees with \cite[Corollary 3.2]{Ha}.
Related to this, it is pointed out in \cite{S88} that the GNS representation of $\phi_\upsilon$ is unitarily equivalent to
the regular representation of $G$ if $|\upsilon| < \sqrt{r/(1-r)}$. % $-\sqrt{r/(1-r)} < \upsilon < \sqrt{r/(1-r)}$.
Since the spectral measures in Theorem~\ref{main} are equivalent for $|\upsilon| \leq \sqrt{r/(1-r)}$,
this remains true even at the critical values $\upsilon = \pm \sqrt{r/(1-r)}$. 

To understand these more clearly, % As a final remark,
we first give a simple proof of the fact due to Haagerup (see \cite{Sz} for a result generalizing this): 

\begin{Proposition}
The average map $E:\C G \to \sA$ is extended to a conditional expectation (also denoted by $E$) $C^*(G) \to A$. 
\end{Proposition}

\begin{proof}
Let $f \in  \C G$ and realize $E(f) \in \sA \subset A$ as a continuous function of $t \in [-1,1]$. For $t \in \sigma_r$, 
\[
  |E(f)(t)| \leq \| E(f)\|' \leq \| f\|' \leq \| f\|,  
\]
where $\|f\|'$ denotes the reduced C*-norm of $f$. 
For $t \in [-1,1]\setminus \sigma_r$, we can find a real $\upsilon$ in the range 
$\sqrt{r/(1-r)} < |\upsilon| \leq 1$ satisfying $t = c_r(\upsilon)$.
In the spectral representation of $\varphi_\upsilon$,
magnify the mass at the point $c_r(\upsilon)$ so that $x^*\varphi_vx \in A^*$ 
represents the Dirac measure at $c_r(\upsilon)$, 
where $x \in A$ is supported by a neighborhood of $c_r(\upsilon)$ as a continuous function on $[-1,1]$.
Note that $\varphi_\upsilon(xE(f)x^*) = \varphi_\upsilon(xx^*E(f)) = E(f)(t)$ particularly. 

Now, given $\epsilon>0$, approximate $x$ by $a \in \sA$ so that $\| xx^* - aa^* \| \leq \epsilon$. 
Then we have $|E(f)(t) - a^*\phi_\upsilon a(E(f))| \leq \epsilon \| E(f)\|$. 
As a positive functional,
$a^*\phi_\upsilon a \in C(G)^*$ satisfies
$|1 - (a^*\phi_\upsilon a)(h_0)| = |\phi_\upsilon(xx^*) - \phi_\upsilon(aa^*)| \leq \epsilon$,
i.e., $1 - \epsilon \leq \| a^*\phi_\upsilon a\| \leq 1 + \epsilon$. 
Putting these together, 
\begin{align*}
 |E(f)(t)| &\leq 
 \epsilon \| E(f)\| + |\phi_\upsilon(aE(f)a^*)| =
             \epsilon \| E(f)\| + |\phi_\upsilon(E(afa^*))|\\
  &= \epsilon \| E(f)\| + |\phi_\upsilon(afa^*)| \leq \epsilon \| E(f)\| + (1+\epsilon) \| f\|.
\end{align*}
% \begin{align*}
%   | E(f)(t)| &\leq |E(f)(t) - h(aE(f)a^*)| + |h(aE(f)a^*)|\\ 
%   &= |E(f)(t) - h(aE(f)a^*)| + |h(E(afa^*))|\\
%   &= |E(f)(t) - h(aE(f)a^*)| + |h(afa^*))|\\
%              &= |E(f)(t) - h(aE(f)a^*)| + |(a^*ha(f))|\\
%              &\leq |E(f)(t) - h(aE(f)a^*)| + \|a^*ha\|\,\|f\|\\
%   &\leq \epsilon + (1+\epsilon) \| f\|. 
% \end{align*}
Since $\epsilon>0$ is arbitrary, this implies $|E(f)(t)| \leq \| f\|$. 
Thus $\| E(f)\| = \sup \{ |E(f)(t)|; -1 \leq t \leq 1\} \leq \| f\|$ and $E$ is extended to
a norm-one projection of $C^*(G)$ onto $A \subset C^*(G)$. 
\end{proof}

As an operator algebraic conditional expectation, $E$ preserves positivity and a positive functional $\varphi$ on $A$
gives rise to a positive functional $\varphi\circ E$ on $C^*(G)$.

Choosing the evaluation functional $\delta_c$ of $A$
at a spectral point $c \in [-1,1]$ as a $\varphi$,
we have a state $\epsilon_c \equiv \delta_c\circ E$ of $C^*(G)$ with the accompanied positive definite function 
% Since states of $C^*(G)$ are rephrased by normalized positive definite functions on $G$,
% $\epsilon_c$ is also described by the associated positive definite function
\[
  \epsilon_c(g) = \delta_c(h_{|g|}) = p_{|g|}(c)\quad (g \in G), 
\]
which is called a spherical function of $G$ and the associated GNS-representation of $G$ is referred to as
a spherical representation of spectral parameter $c$.
Spherical representations are known to be irreducible and mutually inequivalent for different spectral parameters (\cite{FP}). 

Returning to the general positive functional $\varphi$ on $A$,
the GNS-vector and the GNS-representation space of $\varphi\circ E$ are denoted by $(\varphi\circ E)^{1/2}$ and $\sR_\varphi$ 
respectively. The GNS-representation is then called a radial representation and expressed in a form of $C^*(G)$-module on $\sR_\varphi$.
For a spherical representation of spectral parameter $c$, $\sR_{\delta_c}$ is simply denoted by $\sR_c$.

We now observe that $\sR_\varphi$ is a $C^*(G)$-$A$ bimodule (called a radial bimodule) in a natural fashion:
On the level of a dense subspace $\C G (\varphi\circ E)^{1/2}$, the right action of $a \in A$ is realized by
$(x(\varphi\circ E)^{1/2})a = xa (\varphi\circ E)^{1/2}$ ($x \in C^*(G)$), which
in fact gives a well-defined right action of $A$ on $\sR_\varphi$ due to the estimate 
\begin{align*}
  \| xa (\varphi\circ E)^{1/2}\|^2
  &= \varphi(E(a^*x^*xa)) = \varphi(a^*E(x^*x)a) = \varphi(aa^* E(x^*x))\\
  &\leq \|a\|^2 \varphi(E(x^*x)) = \| a\|^2 \| x (\varphi\circ E)^{1/2}\|^2. 
\end{align*}

Given two positive functionals $\varphi$, $\psi$ of $A$,
we use a categorical notation $\Hom(\sR_\varphi,\sR_\psi)$ to express the space of intertwiners between
bimodules $\sR_\varphi$ and $\sR_\psi$. In fact, with these intertwiner spaces as hom-sets, radial representations $\sR_\varphi$ for
various $\varphi$ constitute a W*-category $\sR$ (see \cite{GLR} for W*-categories and \cite{Y07} for further information). 

The irreducibility and inequivalence of spherical representations are then simply expressed by
$\Hom(\sR_s,\sR_t) = \C \delta_{s,t}$ ($s,t \in [-1,1]$). 

Recall that $\tau = \phi_0 = \varphi_0\circ E$ for the choice $v = 0$ and
Cartier's Plancherel formula is a disintegration of $\sR_{\varphi_0}$ into spherical representations.
We shall here generalize this to the radial representation $\sR_\varphi$. 

To this end, we first introduce a continuous field structure on the family $(\sR_t)_{-1 \leq t \leq 1}$ of Hilbert spaces 
by a countable family $(g\epsilon_t^{1/2})$ of sections indexed by $g \in G$: For each $g,g' \in G$, 
$(g\epsilon_t^{1/2}|g'\epsilon_t^{1/2}) = \epsilon_t(g^{-1}g') = p_{|g^{-1}g'|}(t)$ is continuous in $t \in [-1,1]$.

As a consequence it induces a Borel field structure on $(\sR_t)$ and,
under the identification of a positive functional $\varphi$ of $A = C([-1,1])$ with the associated Radon measure $\varphi(dt)$ on $[-1,1]$, 
we can talk about the direct integral of $(\sR_t)$ with respect to $\varphi(dt)$, 
which is denoted\footnote{The square root symbol is used in accordance with von Neumann's suggestion.} by 
\[
  \oint_{[-1,1]} \sR_t\, \sqrt{\varphi(dt)}. 
\]

\begin{Proposition}
The bimodule $\sR_\varphi$ is isometrically isomorphic to $\oint \sR_t\, \sqrt{\varphi(dt)}$ by the correspondence
\[
  x (\varphi\circ E)^{1/2}
  \mapsto
  \oint_{[-1,1]} x \epsilon_t^{1/2}\, \sqrt{\varphi(dt)}
  \quad
  (x \in C^*(G)). 
\]
\end{Proposition}

\begin{proof}
Thanks to the identity 
\[
  \| x(\varphi\circ E)^{1/2}\|^2
  = \varphi(E(x^*x)) = \int_{[-1,1]} \delta_t(E(x^*x))\, \varphi(dt)
  = \int_{[-1,1]} \| x\epsilon_t^{1/2}\|^2\, \varphi(dt),  
\]
the correspondence is isometrically extended to $\sR_\varphi$.

This isometric extension is then surjective. 
In fact, if
\[
  \oint_{[-1,1]} \xi(t)\, \sqrt{\varphi(dt)}
  \in
  \oint_{[-1,1]} \sR_t\, \sqrt{\varphi(dt)}
\]
is orthogonal to
\[
  \oint_{[-1,1]} ga \epsilon_t^{1/2}\, \sqrt{\varphi(dt)}
\]
for every $g \in G \subset C^*(G)$ and $a \in A$, then
\[
  \int_{[-1,1]} (ga\epsilon_t^{1/2}| \xi(t))\, \varphi(dt) = \int_{[-1,1]} \overline{a(t)} (g\epsilon_t^{1/2}| \xi(t))\, \varphi(dt) = 0. 
\]
Since $(g\epsilon_t^{1/2} | \xi(t))$ is integrable with respect to $\varphi(dt)$ and $C([-1,1])$ is weak* dense in $L^\infty(\varphi)$,
$(g\epsilon_t^{1/2}|\xi(t)) = 0$ at $\varphi$-a.e.~$t \in [-1,1]$ for each $g \in G$.
Since $G$ is countable, we can find a $\varphi$-null set $N \subset [-1,1]$ such that
$(g\epsilon_t^{1/2} | \xi(t))= 0$ if $g \in G$ and $t \not\in N$.
Since $G\epsilon_t^{1/2}$ is total in $\sR_t$, this implies $\xi(t) = 0$ for $t \in [-1,1] \setminus N$, i.e.,
$\oint \xi(t)\, \sqrt{\varphi(dt)} = 0$. 
\end{proof}

If we take the positive functional $\varphi_\upsilon$ ($-1 \leq \upsilon \leq 1$) of geometric series, 
the Radon measure $\varphi_\upsilon(dt)$
is explicitly described by Theorem~\ref{main} and the positive definite function $\phi_\upsilon = \varphi_\upsilon\circ E$ of Haagerup
is disintegrated as 
\[
 (\phi_\upsilon)^{1/2} = \oint_{[-1,1]} \epsilon_t^{1/2}\, \sqrt{\varphi_\upsilon(dt)},  
\]
which is reduced to the Plancherel formula in \cite{Ca} for the choice $\upsilon = 0$. 

% In terms of the radial bimodule, the irreducibility of spherical representations is now rephrased as follows:
As a final remark, we give a categorical equivalence between the radial category and the W*-category of GNS-representations of $A$.

For a positive functional $\varphi$ of $A$, let $L^2(\varphi)$ be the GNS-representation of $\varphi$ as an $A$-module and 
$\Hom(L^2(\varphi),L^2(\psi))$ with $\psi$ another positive functional of $A$ be the space of intertwiners, which
constitute a W*-category. 

\begin{Proposition}
  For positive functionals $\varphi$ and $\psi$ of $A$,
  $\Hom(\sR_\varphi,\sR_\psi)$ is naturally isomorphic to $\Hom(L^2(\varphi),L^2(\psi))$.
  In other words, the category of radial bimodules and that of GNS-representations of $A$ are equivalent as W*-categories.   
% The von Neumann algebra $\text{End}({}_{C^*(G)} \sR_A)$ of intertwiners is generated by the right multiplication of $A$ on $\sR$. 
\end{Proposition}

\begin{proof}
  For any functional $\omega \in A^*_+$ which is equivalent to $\varphi+\psi$ as a measure,
  $\oint \sR_t\, \sqrt{\varphi(dt)}$ is isometrically embedded into
  $\oint \sR_t\, \sqrt{\omega(dt)}$ as a bimodule by
  \[
    \oint_{[-1,1]} x\epsilon_t^{1/2}\, \sqrt{\varphi(dt)} \mapsto
    \oint_{[-1,1]} x\epsilon_t^{1/2}\sqrt{\frac{\varphi(dt)}{\omega(dt)}}\, \sqrt{\omega(dt)}
  \]
  and similarly for $\oint \sR_t\, \sqrt{\psi(dt)}$,
  where $\frac{\varphi(dt)}{\omega(dt)}$ denotes the Radon-Nikodym density function.
  % Let $e$ and $f$ be projections to these range subbimodules of $\sR_\omega$.
  % Then $e, f \in \End(\sR_\omega)$ and
  % \[
  %   \Hom(\sR_\varphi,\sR_\psi) \cong \{ T \in \End(\sR_\omega); T = fTe\}. 
  % \]
  % Thus the existence of natural isomorphism is reduced to showing $\End(\sR_\omega) \cong \End(L^2(\omega)$ in view of
  % $\End(L^2(\omega) = L^\infty(\omega)$. 

  The existence of natural isomorphism is then reduced to the case $\varphi = \psi$.
  By the spherical disintegration of $\sR_\varphi$,
  each bounded operator $T$ on $\sR_\varphi$ commuting with the right action of $A$ is decomposed
  as\footnote{The notation $\varphi(dt)^0$ indicates that it depends only on the equivalence class of the measure $\varphi(dt)$.}
  \[
    T = \oint_{[-1,1]} T_t\, \varphi(dt)^0\quad
    \text{with}
    \ 
    T_t \in \cB(\sR_t). 
  \]
  
  Let $\pi_t$ be the spherical representation of $C^*(G)$ on $\sR_t$.
  If $T$ further commutes with countably many decomposable operators
  \[
    \oint_{[-1,1]} \pi_t(g) \, \varphi(dt)^0
  \]
  ($g \in G$), then
  $T_t \in \pi_t(G)'$ for almost any $t \in [-1,1]$.
  Since $\pi_t$ is irreducible, this means that $T_t$ are scalar operators
  for almost any $t$ and $T$ belongs to $L^\infty(\varphi) = \Hom(L^2(\varphi),L^2(\varphi))$. 
\end{proof}

% \textbf{The category of radial representations and the category of GNS-representations of $A$ are equivalent.}

 % Let $\varphi_\upsilon$ ($-1 \leq \upsilon \leq 1$)
 % be the state of geometric series on $A \cong C([-1,1])$ so that % Haagerup's positive definite function is given by
 % $\phi_\upsilon = \varphi_\upsilon\circ E$ corresponds to the positive definite function of Haagerup. %as a state on $C^*(G)$. 
 %   As witnessed in Theorem~\ref{main}, the spectral measure $\mu_\upsilon$ of $\varphi_\upsilon$ then has the following properties. 
 %  \begin{enumerate}
 %  \item
 %    If $|\upsilon| \leq \sqrt{r/(1-r)}$, $\mu_\upsilon$ is equivalent to the Lebesgue measure on $\sigma_r$.
 %  \item
 %    If $\sqrt{r/(1-r)} < |\upsilon| < 1$,
 %    there appears an atomic measure at $c_r(\upsilon) \in (-1,1)\setminus \sigma_r$ adding to a continuous measure equivalent to
 %    the Lebesgue measure on $\sigma_r$.
 %  \item
 %    For $\upsilon = \pm 1$, $\mu_\upsilon$ is the Dirac measure at $\pm 1$. 
 %  \end{enumerate}
 %  Notice that the continuous part of $\mu_\upsilon$ is equivalent to $\mu_0$. 

 %  In view of $\tau = \phi_0$, this is combined with above observation to get the following. 
  % , which supplements Theorem~\ref{SS} (i) and \cite{S88}. 
 % Notice that the critical case $u = \pm \sqrt{r/(1-r)}$ is included here. 

 \begin{Corollary}[\cite{S88,Sz}]\label{unitaryequivalence}~ 
   % Since the tracial state $\tau$ has the same property with the case (i), we have the following alternatives:
   Let $\phi_\upsilon(g) = \upsilon^{|g|}$ for $-1 \leq \upsilon \leq 1$ be a positive definite function of Haagerup. 
  \begin{enumerate}
      \item
      The GNS representation of $\phi_\upsilon$ is unitarily equivalent to $\lambda$ if $|\upsilon| \leq \sqrt{r/(1-r)}$. 
    \item
      The GNS representation of $\phi_\upsilon$ is unitarily equivalent to
      a direct sum of $\lambda$ and the spherical representation for $c = c_r(\upsilon)$ if $\sqrt{r/(1-r)} < |\upsilon| < 1$.
    \item      
    The GNS representation of $\phi_{\pm 1}$ is a character (one-dimensional representation). 
  \end{enumerate}
\end{Corollary}

{\small
  \begin{Remark}~ 
    \begin{enumerate}
    \item One may find some resemblance with the imprimitivity theorem of Rieffel but
      the present isomorphism (or equivalence) is based on the irreducibility of spherical representations and not a consequence
      of a general principle such as systems of imprimitivity. 
      \item The critical case $\upsilon = \pm \sqrt{r/(1-r)}$ in the corollary would be new. 
      \item
        Our original motivation for the present work was to obtain the transition probability in the sense of \cite{Y08}
        between positive definite functions of Haagerup. With the framework in \cite{Y92}, it is expressed by an elliptic integral.
\end{enumerate}
\end{Remark}}

\appendix
\section{Off-Spectrum Region}

Under the correspondence 
\[
  z \mapsto w = \frac{r}{z} + (1-r) z, % - \frac{1}{2} \left( z + \frac{r(1-r)}{z}\right)
\]
the circle $|z| = \sqrt{r/(1-r)}$ is mapped onto the line segment
\[
  \sigma_r = [-2\sqrt{r(1-r)},2\sqrt{r(1-r)}]
\]
with both of open regions $|z| > \sqrt{r/(1-r)}$ and $|z| < \sqrt{r/(1-r)}$ (including $\infty$)
mapped biholomorphically onto $\C \setminus \sigma_r$. 
As their inverses, $\C \setminus \sigma_r$ is mapped biholomorphically onto $|z| > \sqrt{r/(1-r)}$ and
$|z| < \sqrt{r/(1-r)}$ respectively by correspondences 
\[
  z = \frac{w + \sqrt{w^2 - 4r(1-r)}}{2(1-r)}, \quad
  z = \frac{w - \sqrt{w^2 - 4r(1-r)}}{2(1-r)}. 
\]

In these regions, $\pm(2\sqrt{r(1-r)},1]$ for $w$ correspond to $\pm[\frac{r}{1-r},\sqrt{\frac{1-r}{r}})$
for $z$, whence $\C \setminus [-1,1]$ is biholomorphically mapped onto
\[
  D_r =  \Bigl\{ z; |z| > \sqrt{\frac{r}{1-r}}\Bigr\} \setminus
  \Bigl( (\sqrt{r/(1-r)},1]
  \bigcup [-1,-\sqrt{r/(1-r)/r}) \Bigr)
\]
and 
\[
  \Bigl\{ z; |z| < \sqrt{\frac{r}{1-r}} \Bigr\} \setminus
  \Bigl( [(1-r)/r,\sqrt{(1-r)/r}) \bigcup (-\sqrt{r/(1-r)}),-r/(1-r)] \Bigr)
\]
respectively.

Consequently, the extendability of $S(w)$ to $\C \setminus [-1,1]$ in Theorem~2.1 is equivalent to
the extendability of $(\phi(1/z) - r\phi_0)/(r^2/z - (1-r)^2z)$ to $D_r$

\section{Characterization of Stieltjes Transform}
We shall here give a characterization of Stieltjes transforms of compactly supported complex measures.
Let $\mu$ be a complex measure on $\R$ with its support $[\mu]$ compact, 
Then $\mu$ admits moments $(\mu_n)_{n \geq 0}$ and
the Stieltjes transform $S(z)$ ($z \in \C \setminus \R$) of $\mu$ satisfies
\[
  S(z) = \int \frac{1}{t - z}\, \mu(dt) = - \sum_{n \geq 0} \frac{\mu_n}{z^{n+1}}
\]
with the right hand side absolutely convergent for a sufficiently large $|z|$.
Thus $S(z)$ is analytically extended to $\overline{\C} \setminus [\mu]$ so that $S(\infty) = 0$.
Here $\overline{\C} = \C \sqcup \{\infty\}$ denotes the Riemann sphere. 

The complex measure $\mu$ is known to be recovered from $S$ by the inversion formula
\[
  2\pi i \mu(dx) = \lim_{y \to +0} \Bigl( S(x+iy) - S(x-iy) \Bigr)\, dx,
\]
i.e., for any $f \in C_c(\R)$ (the set of compactly supported continuous functions on $\R$), we have 
\[
  2\pi i \int_\R f(x)\, \mu(dx) = \lim_{y \to +0} \int_{-\infty}^\infty f(x) \Bigl( S(x+iy) - S(x-iy) \Bigr)\, dx.
\]
Notice here that the support $[\mu]$ of $\mu$ is easily read off from $S$:
A real point $c$ does not belong to $[\mu]$ if and only if $S$ is holomorphically extended to a neighborhood of $c$ in $\C$.

Based on this fact, the (singularity) support $[\varphi]$ of a holomorphic function $\varphi$ on $\C \setminus \R$ is defined to be
% the complement of
\[
  \{ c \in \overline{\C}; \text{$\varphi$ has no holomorphic extension to a neighborhood of $c$ in $\overline{\C}$} \}. 
\]
Thus $[S]$ is the closure of $[\mu]$ in $\overline{\C}$ if $S$ is the Stieltjes transform of a complex measure $\mu$.

\begin{Theorem}\label{characterization}
  A holomorphic function $S(z)$ on $\C \setminus \R$ is the Stieltjes transform of
  a compactly supported complex measure on $\R$ if and only if the following conditions are satisfied.
%  satisfying the three conditions is the Stieltjes transform of the measure $\mu$ in (iii). 
  \begin{enumerate}
      \item
  $S(z)$ is analytically extended to a neighborhood of $\infty$ in the Riemann sphere $\overline{\C} = \C \sqcup \{ \infty\}$ so that
  $S(\infty) = 0$.
\item
  The limit
  \[
   \lim_{y \to +0} \frac{1}{2\pi i} \Bigl( S(x+iy) - S(x-iy)\Bigr)\,dx
  \]
  exists as a complex Radon measure on $C_c(\R)$. 
\end{enumerate}
\end{Theorem}

\begin{proof}
  The necessity is already checked above. 
%  The conditions are necessary as seen above.

  To see the sufficiency, let $S_\mu$ be the Stieltjes transform of the complex measure $\mu$ in (ii).
  By the condition (i), the support $[\mu]$ of $\mu$ is compact and therefore 
  $S_\mu$ satisfies the properties in (i) and (ii) as well. Thus their difference $\varphi = S - S_\mu$
  also meets the conditions in such a way that
  \[
\lim_{y \to 0} \int_\R f(t) \Bigl(\varphi(t+iy) - \varphi(t-iy)\Bigr)\, dt = 0
  \]
  for any $f \in C_c(\R)$ and the problem is to show that $\varphi = 0$.

  To this end, we simply modify the argument in \cite[\S 8]{Ru} as follows:
  Regarding $\varphi(z) = \varphi(x+iy)$ as a function of $x \in \R$ with $0 \not= y \in \R$ a parameter,
  take the moving average with $f \in C_c(\R)$ to obtain a holomorphic function 
  \[
   % (f*\varphi)(z) = (f*\varphi)(x+iy) =
    \varphi_f(z) = \int_{-\infty}^\infty f(t) \varphi(x+t + iy)\, dt
  \]
  of $z = x + iy \in \C \setminus \R$, which is analytically extended to  $\overline{\C} \setminus ([\varphi] - [f])$
  and vanishes at $\infty \in \overline{\C}$ from the condition (i). Here $[f]$ denotes the support of $f$. 
  
  Now focus on the harmonic function $\varphi_f(z) - \varphi_f(\overline{z})$ of $z \in \C \setminus ([\varphi] - [f])$.
  We shall show that this is continuously extended to $\C$ in such a way that it vanishes on $\R$.

  Choose $\delta>0$ and a bounded interval $[a,b]$ in $\R$ large enough so that $[f] \subset [-\delta,\delta]$ and
  $[\varphi] \subset [a+3\delta,b-3\delta]$.
  In what follows, we identify the subspace $\{ g \in C_c(\R); [g] \subset [a,b]\}$ of $C_c(\R)$ with
  the set $C_0(a,b)$ consisting of continuous functions on $[a,b]$ vanishing at the boundary $\{a,b\}$ 
  in an obvious manner and furnish these with the topology by the uniform norm.

  Define a bounded linear functional $\Phi_y$ ($0 \not= y \in \R$) on $C_0(a,b)$ by
  \[
    \Phi_y(g) = \int_a^b g(t) (\varphi(t+iy) - \varphi(t-iy))\, dt 
  \]
  for $g \in C_0(a,b)$, which is weak*-continuous in $0 \not= y \in \R$ and satisfies
  \[
    \lim_{y \to 0} \Phi_y(g) = 0 = \lim_{y \to \pm\infty} \Phi_y(g). 
  \]
  
  In view of $[\varphi] - [f] \subset [a+2\delta,b-2\delta]$, $\varphi_f$ is holomorphic on $\overline{\C} \setminus [a+2\delta,b-2\delta]$,
  whereas, if $(T_xf)(t) = f(t-x)$ denotes the translate of $f$ by $x \in \R$, $T_xf \in C_0(a,b)$ for $x \in [a+\delta,b-\delta]$ with
  $T_xf$ norm-continuous in $x \in [a+\delta,b-\delta]$ and we have 
  \[
    \varphi_f(z) - \varphi_f(\overline{z}) = \Phi_y(T_xf)
    \quad
    (x \in [a+\delta,b-\delta], 0 \not= y \in \R).
  \]
  
  % which is continuously extended to
  % expressed as a combination of two operations. Given $0 \not= y \in \R$,
  % introduce a linear functional $\Phi_y$ of $f \in C_c(\R)$
  % by
  % \[
  %   \Phi_y(f) = \int_{-\infty}^\infty f(t) \Bigl( \varphi(-t+iy) - \varphi(-t-iy) \Bigr)\, dt, 
  % \]
  % which is well-defined and converges weakly to $0$ as $y \to 0$ because $\varphi(z)$ vanishes at $\infty$,
  % and define the (right) translation $T_xf$ of $f$ under $x \in \R$ by $(T_xf)(t) = f(t+x)$. Then we have 
  % \[
  %   (f*\varphi)(x+iy) - (f*\varphi)(x-iy) = \Phi_y(T_xf).
  % \]
  
  By the Banach-Steinhaus theorem, the functional norm of $\Phi_y$ is uniformly bounded in $0 \not= y \in \R$,
  which is combined with uniform norm continuity of $T_xf$ in $x \in [a+\delta,b-\delta]$ to see that
  $\Phi_y(T_xf)$ is extended to a continuous function of $(x,y) \in [a+\delta,b-\delta]\times \R$ so that $\Phi_0(T_xf) = 0$.
  Since $\varphi_f(z) - \varphi_f(\overline{z})$ is continuous in $z \in \C \setminus [a+2\delta,b-2\delta]$,
  we have proved that the harmonic function $\varphi_f(z) - \varphi_f(\overline{z})$ on $\C \setminus \R$
  is continuously extended to $\overline{\C}$ with
  the continuous extension vanishing on $\R \sqcup \{ \infty\}$. 

  Now the maximum principle for harmonic functions is applied on the upper or lower half plane to conclude that
  $\varphi_f(z) - \varphi_f(\overline{z})$ vanishes identically.
  
  % Since the function $\Phi_y(T_xf)$ is harmonic on the upper half plane $y > 0$ with its boundary function vanishing including
  % the infinite point, the maximal principle for harmonic functions implies that $\Phi_y(T_xf) = 0$ ($x \in \R$, $y > 0$).
  % Likewise $\Phi_y(T_xf)$ vanishes identically on the lower halfplane $y < 0$.

  Consequently we have $\varphi_f(z) = \varphi_f(\overline{z})$ for any $f \in C_c(\R)$. 
  % \[
  %   (f*\varphi)(z) = (f*\varphi)(\overline{z}).
  % \]
  % for any $f \in C_0(\R)$.
  Choosing a sequence $(f_n)$ of approximate delta functions in $C_c(\R)$ and taking the limit $n \to \infty$,
  we obtain $\varphi(z) = \varphi(\overline{z})$ ($z \in \C \setminus \R$) and therefore $\varphi(z)$ is constant
  with its value $\varphi(\infty) = 0$. 
\end{proof}

\begin{Remark}
  It is standard to identify complex measures on $\R$ with bounded linear functionals on the C*-algebra $C_0(\R)$ of
  continuous functions vanishing at infinity (Riesz representation theorem). By a simple modification,
  we can also identify locally complex measures on $\R$ with locally bounded linear functionals of $C_c(\R)$, which is used
  in the condition (ii) and gives rise to compactly supported complex measures by the condition (i).

  By reversing the order of usage, we can even work within the standard identification. In fact, the condition (i) implies
  $S(x+iy) - S(x-iy) = O(1/x^2)$ for each $y > 0$ and $S(x+iy) - S(x-iy)$ ($y > 0$) are integrable functions.
  Thus $(S(x+iy) - S(x-iy))dx$ ($y>0$) define complex measures on $\R$ and the condition (ii) is then strengthened
  to the weak convergence on $C_0(\R)$. % this time. 
\end{Remark}

\end{document}